\documentclass[a4paper, 12pt, reqno]{amsart}
\usepackage[utf8]{inputenc}
\usepackage{amssymb,mathtools,cite,enumerate,color,eqnarray,hyperref,amsfonts,amsmath,amsthm,setspace,tikz,verbatim,charter,booktabs,multirow,xcolor,latexsym}
\usepackage[a4paper,top=3cm,bottom=3cm,left=2.5cm,right=2.5cm]{geometry}
\usepackage{graphics,graphicx}
\usepackage{txfonts}
\usepackage{youngtab}
\usepackage{tabularx}
\usepackage{ytableau}
\usetikzlibrary{tikzmark,shapes.geometric, arrows.meta, positioning, fit,decorations.pathreplacing}
\usepackage{subcaption}
\usepackage[none]{hyphenat}
\usepackage[skip=8pt plus 1pt, indent=15pt]{parskip}

\numberwithin{equation}{section}
\definecolor{ao(english)}{rgb}{0.0, 0.0, 0.6}
\hypersetup{colorlinks=true, linkcolor=ao(english),citecolor=ao(english)}

\usepackage[normalem]{ulem}

\theoremstyle{plain}
\newtheorem{theorem}{Theorem}
\numberwithin{theorem}{section}

\newtheorem*{corollary*}{Corollary}
\newtheorem*{Example*}{Example}

\newtheorem{conjecture}[theorem]{Conjecture}
\theoremstyle{definition}

\newtheorem*{def*}{Definition}
\newtheorem*{theorem*}{Theorem}

\newtheorem*{definition*}{Definition}

\theoremstyle{remark}

\allowdisplaybreaks

	\title[Fixed Perimeter]{A Proof of a Conjecture on Fixed Perimeter Partitions}

        \author[P. J. Mahanta]{Pankaj Jyoti Mahanta}
\address{Department of Mathematical Sciences, Tezpur University, Assam, India, PIN-784028}
\email{pjm2099@gmail.com, msp25007@tezu.ac.in}

	
\dedicatory{The author dedicates this paper to Prof.~Nayandeep Deka Baruah, one of his favorite Assamese mathematicians.}
	
\begin{document}

\begin{abstract}
Finding fixed perimeter analogues of various partition theoretic identities and inequalities has recently emerged as an active area of research. Gray, Payne, Swisher, and Watson [\textit{Discrete Math.}, 2026] established several fixed perimeter analogues of partition theoretic results inspired by Euler's celebrated partition identity. Very recently, in a separate work [\textit{ar{X}iv:2608.00421}, 2026], they explored fixed perimeter analogues of inequalities related to parity biases. Introducing the concept of parity bias, B. Kim, E. Kim, and Lovejoy [\textit{Eur. J. Comb.}, 2020] conjectured that $pd_o(n)>pd_e(n)$ for all $n\ge 20$, where $pd_o(n)$ (respectively, $pd_e(n)$) denote the number of partitions of $n$ into distinct parts having more odd parts (respectively, even parts) than even parts (respectively, odd parts). The author, together with Banerjee, Bhattacharjee, Dastidar, and Saikia [\textit{Eur. J. Comb.}, 2022], proved this conjecture. Gray, Payne, Swisher, and Watson conjectured that a fixed perimeter analogue of this inequality holds for all $n\ge 9$. In this paper, we confirm their conjecture.
\end{abstract}
	\maketitle

	\section{Introduction}
\begin{figure}[h]
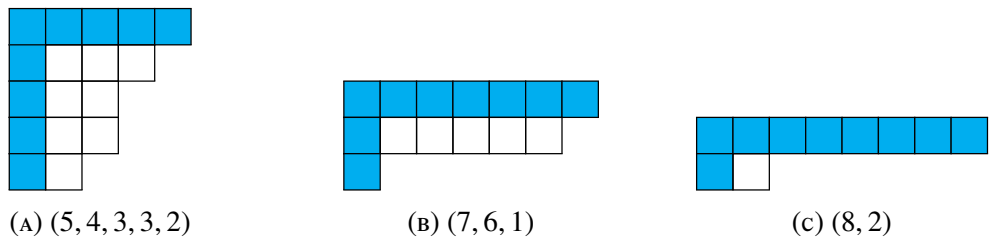

\centering
\ytableausetup{boxsize=1.1em}
\begin{minipage}[b]{0.3\textwidth}
		\centering
\begin{ytableau}
	*(cyan)~ & *(cyan)~ & *(cyan)~& *(cyan)~ & *(cyan)~ \\
	*(cyan)~ & ~ & ~ & ~ \\
	*(cyan)~ & ~ & ~\\
	*(cyan)~ & ~ & ~\\
	*(cyan)~ & ~
\end{ytableau}
\subcaption{$(5,4,3,3,2)$}
\end{minipage}
\begin{minipage}[b]{0.3\textwidth}
		\centering
\begin{ytableau}
	*(cyan)~ & *(cyan)~ & *(cyan)~& *(cyan)~ & *(cyan)~ & *(cyan)~ & *(cyan)~\\
	*(cyan)~ & ~ & ~ & ~ & ~ & ~ \\
	*(cyan)~
\end{ytableau}
\subcaption{$(7,6,1)$}
\end{minipage}
\begin{minipage}[b]{0.3\textwidth}
		\centering
\begin{ytableau}
	*(cyan)~ & *(cyan)~ & *(cyan)~& *(cyan)~ & *(cyan)~ & *(cyan)~ & *(cyan)~ & *(cyan)~\\
	*(cyan)~ & ~
\end{ytableau}
	\subcaption{$(8,2)$}
\end{minipage}
\caption{Three partitions, each with a perimeter of 9. The boxes comprising the largest hook (perimeter) are highlighted.}
\label{fig:perimeter}
\end{figure}

A \textit{partition} of a positive integer $n$ is a sequence of positive integers $\pi=(\pi_1, \pi_2, \ldots, \pi_r)$ such that $\pi_1\geq \pi_2\geq \cdots \geq \pi_r$ and $\sum\limits_{i=1}^r\pi_i = n$. The numbers $\pi_1, \pi_2, \ldots, \pi_r$ are called the \textit{parts} of the partition $\pi$. We denote by $p(n)$ the number of partitions of $n$. For a foundational and extensive background on partition theory, the reader is referred to Andrews' classic book~\cite{andrews1998theory}. The \textit{Young diagram} of a partition $(\pi_1, \pi_2, \ldots, \pi_r)$ is a left-justified array of boxes, where the $i$-th row from the top contains $\pi_i$ boxes. The \textit{hook length} of a box in a Young diagram is the sum of the number of boxes directly to its right, the number of boxes directly below it, and 1 (for the box itself). Representing partitions using Young diagrams and exploring their properties through diagrammatic analysis is a vast sub-field of partition theory. One such property is the \textit{perimeter} of a partition, which is defined to be the largest hook length of the partition. Figure~\ref{fig:perimeter} depicts three partitions, each having a perimeter of 9.

In 2016, Straub~\cite{straub2016core} proved a fixed perimeter analogue of Euler’s famous and remarkably intriguing partition identity,
$$o(n)=d(n), \quad \text{for all } n\ge 1,$$
where $o(n)$ and $d(n)$ denote the number of partitions of $n$ into odd parts and distinct parts, respectively. His work inspired researchers to prove fixed perimeter analogues of various partition identities and inequalities. In 2018, Fu and Tang~\cite{fu2018partitions} proved that the number of partitions with perimeter $n$ is $2^{n-1}$. The connection between fixed perimeter and Fibonacci numbers and Tribonacci numbers has also attracted significant interest (see~\cite{gray2026_arxiv, straub2016core}). Several recent developments on this topic can be found in~\cite{chen-integer, gray24proceedings, lin2022combinatorics, waldron}.

In~\cite{gray2026}, Gray, Payne, Swisher, and Watson proved several fixed perimeter analogues of partition results inspired by Euler’s aforementioned identity. In~\cite{gray2026_arxiv}, they extended their exploration to parity biases, a concept introduced by B. Kim, E. Kim, and Lovejoy~\cite{kim2020parity}. Motivated by B. Kim, E. Kim, and Lovejoy's work, the author, together with Banerjee, Bhattacharjee, Dastidar, and Saikia~\cite{banerjee2022parity}, proved the following bias (originally conjectured by B. Kim, E. Kim, and Lovejoy) along with several other parity biases.

\begin{theorem}\cite[Theorem 1.4]{banerjee2022parity}
		For all $n\ge 20$, we have
	$$pd_o(n)>pd_e(n),$$
	where $pd_o(n)$ (respectively, $pd_e(n)$) denote the number of partitions of $n$ into distinct parts having more odd parts (respectively, even parts) than even parts (respectively, odd parts).
\end{theorem}

Let $rd_o(n)$ and $rd_e(n)$ be the fixed perimeter analogues of $pd_o(n)$ and $pd_e(n)$, respectively. In \cite{gray2026_arxiv}, Gray, Payne, Swisher, and Watson posed the following conjecture.

\begin{conjecture}\cite[Conjecture 1.4]{gray2026_arxiv} \label{conjecture:perimeter_bias_2026}
	For all $n\ge 9$, we have
	$$rd_o(n)>rd_e(n).$$
\end{conjecture}

In this paper, we prove the following theorem.

\begin{theorem}\label{thm:perimeter_bias_2026}
	Conjecture \ref{conjecture:perimeter_bias_2026} is true.
\end{theorem}

Before proceeding to the proof, we recall the following notation:
\begin{itemize}
	\item $RD_o(n)$ and $RD_e(n)$ denote the sets of partitions counted by $rd_o(n)$ and $rd_e(n)$, respectively.
	\item $\omega(\pi)$ denotes the number of odd parts in $\pi$ minus the number of even parts in $\pi$.
\end{itemize}

\section{Proof of Theorem~\ref{thm:perimeter_bias_2026}}

We first divide $RD_e(n)$ into three disjoint subsets based on the following three cases. Then, by applying explicit mappings to each partition $\pi\in RD_e(n)$, we construct a corresponding partition $\pi^\prime\in RD_o(n)$.

\noindent \textbf{Case A: The smallest part is $\ge 3$.}
In this case, for each partition $\pi\in RD_e(n)$, we reduce the size of all parts by 1 and then add a part of size 1. Here, if $\omega(\pi)=-m$ for a positive integer $m$, then we have $\omega(\pi^\prime)=m+1$.

\noindent \textbf{Case B: The smallest part is 1.}
We apply the restriction function of $\varphi_1$ defined in \cite{gray2026_arxiv} to each partition. In this case as well, if $\omega(\pi)=-m$, then $\omega(\pi^\prime)=m+1$. 

After applying \textbf{Case A}, the partitions $\lambda$ with smallest part 1 and $\omega(\lambda)=1$ remain in $RD_o(n)$. Similarly, after applying \textbf{Case B}, the partitions $\lambda$ with smallest part $\ge 3$ and $\omega(\lambda)=1$ remain in $RD_o(n)$. Moreover, the partitions $\lambda$ with smallest part 2 remain in $RD_o(n)$. We denote these disjoint subsets of partitions by $RD_o^1(n)$, $RD_o^2(n)$ and $RD_o^3(n)$, respectively.

\noindent \textbf{Case C: The smallest part is 2.}
We divide these partitions into three disjoint subsets based on the value of $\omega(\pi)$ for each $\pi\in RD_e(n)$

\noindent \textbf{Case C1: $\omega(\pi)=-1$.}
In this case, we reduce the size of the smallest part of each partition by 1. Then the image partitions lie in $RD_o^1(n)$. The unmapped partitions in $RD_o^1(n)$ each contain the part 2. We denote this set of remaining partitions by $RD_o^{1^C}(n)$.

\noindent \textbf{Case C2: $\omega(\pi)=-2$.}
In this case, we apply the restriction function of $\varphi_2$ defined in \cite{gray2026_arxiv} to each partition. Here, the image partitions lie in $RD_o^2(n)$. Each unmapped partition in $RD_o^2(n)$ contains 3 as a part. We denote this set of remaining partitions by $RD_o^{2^C}(n)$.

\noindent \textbf{Case C3: $\omega(\pi)\le -3$.} Here, we divide the partitions into three subcases based on the size of the two largest parts of $\pi=(\pi_1, \pi_2, \ldots, \pi_r)$ in $RD_e(n)$. Prior to explaining these subcases, we divide the set $RD_o^3(n)$ into the following three disjoint subsets:
\begin{itemize}
	\item $RD_o^{3a}(n) := \{\lambda\in RD_o^3(n) \mid 3 \text{ is a part}\}$,
	\item $RD_o^{3b}(n) := \{\lambda\in RD_o^3(n)\setminus RD_o^{3a}(n) \mid \text{ the largest part is odd}\}$,
	\item $RD_o^{3c}(n) := \{\lambda\in RD_o^3(n)\setminus RD_o^{3a}(n) \mid \text{ the largest part is even}\}$.
\end{itemize}

\noindent\textbf{Case C3(i): $\pi_1-\pi_2=1$.}
Here, we first remove the part $\pi_1$. Next, we increase each part by 1, excluding the parts of size 2, and add a part of size 3. The resulting image partitions lie in $RD_o^{3a}(n)$, satisfying $\omega(\pi^\prime)\ge 1$.

\noindent\textbf{Case C3(ii): $\pi_1-\pi_2\ge 2$ and $\pi_1$ is odd.}
Here, we increase each part by 1, with the exception of the largest part and the part of size 2. The resulting images lie in $RD_o^{3b}(n)$.

\noindent\textbf{Case C3(iii): $\pi_1-\pi_2\ge 2$ and $\pi_1$ is even.}
We again divide this subcase into the following three subcases. Here, we classify the partitions based on the value of $\omega(\pi)$.

\noindent\textbf{Case iii-a: $\omega(\pi)\le -5$.} Here, we apply the same function as defined in \textbf{Case C3(ii)}. Then, the resulting images lie in $RD_o^{3c}(n)$.

\noindent\textbf{Case iii-b: $\omega(\pi)=-3$.} Here, we increase each part by $1$, with the exception of the largest part. The resulting images lie in $RD_o^{2^C}(n)$.

\noindent\textbf{Case iii-c: $\omega(\pi)=-4$.}
Let $\pi=(\pi_1, \pi_2, \ldots, \pi_r)\in RD_e(n)$ belong to this sub-case. Because $\omega(\pi)=-4$, the number of parts $r$ must be even. We first reduce $\pi_1$ by 1 and add a part of size 1. Since $\pi_r=2$ in this case, the resulting partition is of the form
$$(\pi_1-1, \pi_2, \ldots, \pi_{r-1},2,1),$$
and we have $\omega((\pi_1-1, \pi_2, \ldots, \pi_{r-1},2,1))=-1$. Therefore, among the parts $\pi_2,\pi_3, \ldots, \pi_{r-1}$, exactly $\dfrac{r}{2}$ parts are even and $\dfrac{r-4}{2}$ parts are odd. To obtain a partition $\pi^\prime$ satisfying $\omega(\pi^\prime)=1$ from $(\pi_1-1, \pi_2, \ldots, \pi_{r-1},2,1)$, we need to transform one of the even parts in $\{\pi_2,\pi_3, \ldots, \pi_{r-1}\}$ into an odd part. That is, we need to make number of even parts and odd parts equal, namely $\dfrac{r-2}{2}$ each. Note that the resulting partition $\pi^\prime$ belongs to $RD_o^{1^C}(n)$.

Now, since $\pi_1, \pi_2, \ldots, \pi_r$ are distinct, we have $\pi_{r-1}\ge 3$ and $\pi_2\ge r$. An elementary calculation shows that the number of ways to choose $\dfrac{r}{2}$ even numbers and $\dfrac{r-4}{2}$ odd numbers from $\{3,4,\ldots,s\}$, where $s\ge r$ is a positive integer, is less than the number of ways to choose $\dfrac{r-2}{2}$ even numbers and $\dfrac{r-2}{2}$ odd numbers from the same set.

Routine verification for small values of $n$ completes the proof.

For $n=13$, \textbf{Case C} is presented in Table \ref{table-13}. For this fixed perimeter, there are no pre-images in Cases \textbf{C3(ii)}, \textbf{iii-a}, and \textbf{iii-b}. Therefore, these sub-cases are presented in Table \ref{table-14} for the fixed perimeter $n=14$.

	\begin{table}[h]
	\centering
	\small
	\begin{tabular}{|c|c|c|m{6cm}|m{6cm}|}
		\hline
		\multicolumn{3}{|c|}{\textbf{Cases}} & \textbf{Pre-images} & \textbf{Images} \\ \hline
		\multicolumn{3}{|c|}{\textbf{C1}} & $(11, 10, 2)$, $(11, 8, 2)$, $(11, 6, 2)$, $(11, 4, 2)$, $(9, 8, 7, 6, 2)$, $(9, 8, 7, 4, 2)$, $(9, 8, 6, 5, 2)$, $(9, 8, 6, 3, 2)$, $(9, 8, 5, 4, 2)$, $(9, 8, 4, 3, 2)$, $(9, 7, 6, 4, 2)$, $(9, 6, 5, 4, 2)$, $(9, 6, 4, 3, 2)$ & $(11, 10, 1)$, $(11, 8, 1)$, $(11, 6, 1)$, $(11, 4, 1)$, $(9, 8, 7, 6, 1)$, $(9, 8, 7, 4, 1)$, $(9, 8, 6, 5, 1)$, $(9, 8, 6, 3, 1)$, $(9, 8, 5, 4, 1)$, $(9, 8, 4, 3, 1)$, $(9, 7, 6, 4, 1)$, $(9, 6, 5, 4, 1)$, $(9, 6, 4, 3, 1)$ \\
		\hline
		\multicolumn{3}{|c|}{\textbf{C2}} & $(12, 2)$, $(10, 9, 8, 2)$, $(10, 9, 6, 2)$, $(10, 9, 4, 2)$, $(10, 8, 7, 2)$, $(10, 8, 5, 2)$, $(10, 8, 3, 2)$, $(10, 7, 6, 2)$, $(10, 7, 4, 2)$, $(10, 6, 5, 2)$, $(10, 6, 3, 2)$, $(10, 5, 4, 2)$, $(10, 4, 3, 2)$, $(8, 7, 6, 5, 4, 2)$, $(8, 7, 6, 4, 3, 2)$, $(8, 6, 5, 4, 3, 2)$ & $(13)$, $(11, 10, 9)$, $(11, 10, 7)$, $(11, 10, 5)$, $(11, 9, 8)$, $(11, 9, 6)$, $(11, 9, 4)$, $(11, 8, 7)$, $(11, 8, 5)$, $(11, 7, 6)$, $(11, 7, 4)$, $(11, 6, 5)$, $(11, 5, 4)$, $(9, 8, 7, 6, 5)$, $(9, 8, 7, 5, 4)$, $(9, 7, 6, 5, 4)$ \\
		\hline
		\multirow{5}{*}{\textbf{C3}} & \multicolumn{2}{|c|}{\textbf{C3(i)}} & $(9, 8, 6, 4, 2)$ & $(9, 7, 5, 3, 2)$ \\ \cline{2-5}
		& \multicolumn{2}{|c|}{\textbf{C3(ii)}} & none &  \\ \cline{2-5}
		& \multirow{3}{*}{\textbf{C3(iii)}} & \textbf{iii-a} & none &  \\ \cline{3-5}
		& & \textbf{iii-b} & none &  \\ \cline{3-5}
		& & \textbf{iii-c} & $(10, 8, 6, 2)$, $(10, 8, 4, 2)$, $(10, 6, 4, 2)$ & Partitions in $RD_o^{1^C}(13)$ are: $(9, 4, 3, 2, 1)$, $(9, 6, 3, 2, 1)$, $(9, 8, 3, 2, 1)$, $(9, 5, 4, 2, 1)$, $(9, 7, 4, 2, 1)$, $(9, 6, 5, 2, 1)$, $(9, 8, 5, 2, 1)$, $(9, 7, 6, 2, 1)$, $(9, 8, 7, 2, 1)$ \\
		\hline
	\end{tabular}
	\caption{Partitions and their images for \textbf{Case C} with the fixed perimeter $n=13$.}
	\label{table-13}
\end{table}

\begin{table}[h]
	\centering
	\small
\begin{tabular}{|c|m{7cm}|m{7cm}|}
	\hline
	{\textbf{Cases}} & \textbf{Pre-images} & \textbf{Images}\\
	\hline
	\textbf{C3(ii)} & $(11, 8, 6, 2)$, $(11, 8, 4, 2)$, $(11, 6, 4, 2)$ & $(11, 9, 7, 2)$, $(11, 9, 5, 2)$, $(11, 7, 5, 2)$ \\
	\hline
	\textbf{iii-a} & $(10, 8, 6, 4, 2)$ & $(10, 9, 7, 5, 2)$ \\
	\hline
	\textbf{iii-b} & $(12, 10, 2)$, $(12, 8, 2)$, $(12, 6, 2)$, $(12, 4, 2)$, $(10, 8, 7, 6, 2)$, $(10, 8, 7, 4, 2)$, $(10, 8, 6, 5, 2)$, $(10, 8, 6, 3, 2)$, $(10, 8, 5, 4, 2)$, $(10, 8, 4, 3, 2)$, $(10, 7, 6, 4, 2)$, $(10, 7, 5, 4, 2)$, $(10, 6, 4, 3, 2)$, $(10, 6, 5, 4, 2)$ & $(12, 11, 3)$, $(12, 9, 3)$, $(12, 7, 3)$, $(12, 5, 3)$, $(10, 9, 8, 7, 3)$, $(10, 9, 8, 5, 3)$, $(10, 9, 7, 6, 3)$, $(10, 9, 7, 4, 3)$, $(10, 9, 6, 5, 3)$, $(10, 9, 5, 4, 3)$, $(10, 8, 7, 5, 3)$, $(10, 8, 6, 5, 3)$, $(10, 7, 5, 4, 3)$, $(10, 7, 6, 5, 3)$ \\
	\hline
\end{tabular}
	\caption{Partitions and their images for Cases \textbf{C3(ii)}, \textbf{iii-a}, and \textbf{iii-b} with the fixed perimeter $n=14$.}
\label{table-14}
\end{table}

	\section*{Statements and Declarations}
	
	\textbf{Data Availability:} This manuscript does not contain any associated data.
	
	\textbf{Competing Interests:} The author declares no competing interests.
	
	\textbf{Funding Information:}
    The author was partially supported by an institutional fellowship for doctoral research from Tezpur University, Assam, India.


\end{document}